# Urban Water Sprinkler Routing: A Multi-Depot Mixed Capacitated Arc Routing Problem Incorporating Real-Time Demands

Hongtai Yang[a], Luna Liu[a], Ke Han[b,*], Boyi Lei[a]

[a]School of Transportation and Logistics, National Engineering Laboratory of Integrated Transportation Big Data Application Technology, National United Engineering Laboratory of Integrated and Intelligent Transportation, Institute of System Science and Engineering, Southwest Jiaotong University, Chengdu, China, 611756
[b]School of Economics and Management, Southwest Jiaotong University, Chengdu, China, 611756

*Corresponding author, Email address: kehan@swjtu.edu.cn

## Abstract

Fugitive road dust (FRD), as one of the major pollutants in the city, poses great harm to the environment and the physical health of citizens. A common countermeasure adopted by government agencies is employing on-road water trucks (sprinklers) to spray water (sprinkle) on urban streets to reduce the FRD. Currently, the traveling routes of sprinklers are usually planned based on drivers' experience, which may lead low operation efficiency and could not respond to the real-time sprinkling demands. To address these issues, this study formulates the routes planning of sprinklers as a multi-depot mixed capacitated arc routing problem with real-time demands with the aim of minimizing the sprinklers' travel distance. We develop an improved adaptive large neighborhood search (ALNS) algorithm that incorporates a tabu-list and a perturbation mechanism to solve this problem. Furthermore, a problem-specific acceleration mechanism is designed to reduce unnecessary search domains to improve the efficiency of the algorithm. Empirical experiments are conducted based on various scenarios and the results demonstrate that the proposed algorithm generates solutions that are superior or at least comparable to the solutions generated by the traditional ALNS algorithm but with significantly lower computation time. Sensitivity analysis is conducted to explore the effects of relevant parameters on the results. This study is the first to incorporate real-time FRD pollution information, gathered through multiple data sources via IoT technology, into urban sprinkling operations, extending the traditional CARP from a tactical planning to a real-time operational environment. A real-world implementation case is also presented.

## 1. Introduction

Fugitive road dust (FRD) is one of the major pollutants emitted by motor vehicles (Zhang et al., 2023; Zhu et al., 2009). FRD dominates the traffic-related $PM_{10}$ emissions, constituting an average contribution of 62% to traffic $PM_{10}$ and 22% to traffic $PM_{2.5}$ in some Chinese megacities (Zhang et al., 2020). FRD also serves as a vital carrier of hazardous substances, including heavy metals and polycyclic aromatic hydrocarbons (PAHs), which induce substantial risks to human health (Chen et al., 2019; Khan and Strand, 2018; Li et al., 2020). Therefore, governments agencies around the globe have made significant efforts to reducing FRD.

Water-sprinkling is one of the most widely adopted measures to reduce FRD (Gulia et al., 2019) and is capable of removing approximately 77% to 100% of $PM_{10}$ and 68% to 87% of $PM_{2.5}$ (Zhang et al., 2023). Thus, it is adopted by most cities in China. During sprinkling operations, the water-filled sprinkler usually starts from the depot and follows a designated route. The typical capacity of the water-tanks ranges from 5t to 12t, and the maximum service time of a sprinkler with a full tank of water is about 2 hours. It is crucial to reduce the non-operational (no sprinkling required) travel distances as they are associated with extended working time and extra energy consumption. In addition, considering the large number of on-going construction projects in medium and large cities, the significant amount of FRD associated with earthwork/concrete transportation activities could potentially exacerbate air pollution (Alshetty and M., 2022; Alshetty and Nagendra, 2022). As those transportation activities and associated FRD risks occur in real time and randomly throughout the road network, it is necessary to respond to unforeseen sprinkling demands in a timely manner to minimize environmental and public health hazards. This requires the route planning process to flexibly switch between fixed-route execution (tactical) and random demands response (real-time).

Aiming at improving the efficiency of sprinkling operations (minimize sprinklers' total travel distance) and timely respond to real-time sprinkling demands, this paper formulates the sprinkling route planning as a multi-depot mixed capacitated arc routing problem with real-time demands. To our knowledge, this work is the first to address real-time demands in the routing of sprinklers. The problem is more challenging than the conventional capacitated arc routing problem (CARP) in the following aspects: (1) The demand edges are a mix of directional and non-directional ones; and (2) the routes may be dynamically altered to meet real-time demands. For these reasons, this problem is termed Dynamix-CARP in this paper. To address these challenges, an improved adaptive large neighborhood search (ALNS) heuristic algorithm is developed specifically for the Dynamix-CARP. In the route planning and adjustment stages, we devise specialized destroy and repair operators to efficiently explore solution spaces. To enhance algorithm efficiency, an acceleration mechanism is introduced to reduce search domains during the repair process, and a tabu-list is established to prevent duplicate searches. Furthermore, to expand the search space, the algorithm strategically incorporates perturbations and allows suboptimal solutions to participate in optimization when multiple iterations fail to find improved solutions. Computational experiments across diverse scenarios in a real-world case study demonstrate the efficiency of the proposed algorithm in planning fixed routes (tactical) and making real-time route adjustments (operational), which effectively addresses fixed demands while dynamically responding to emerging demands in realistic urban sprinkling operations.

The remainder of this paper is structured as follows. Section 2 reviews the related literatures. Section 3 presents the formulation of the proposed model. Section 4 introduces the proposed ALNS algorithm. Section 5 presents results of computational experiments and sensitivity analysis. Section 6 concludes the paper.

## 2. Literature review

This section reviews the literatures on the routes planning of sprinklers, variants of CARP, and ALNS algorithm.

Studies on the routes planning of sprinklers can be divided into two main categories: the routes planning on urban streets and the routes planning in specialized places such as open pit mines. Among them, Zhu et al. (2007) focused on the route planning on urban streets and proposed an improved genetic algorithm (GA) to solve the multi-depot capacitated arc routing problem (MDCARP). Yu et al. (2019) explored the split-delivery mixed capacitated arc routing problem (SDMCARP) in the context of sprinkling on urban streets and introduced a forest-based tabu search algorithm to accelerate the solution. Considering the sprinkling operations in open pit mines, Li et al. (2008) designed two heuristics, a minimum cost flow-based heuristic and a set partitioning-based heuristic to solve the MDCARP. Also for the sprinkling in open pit mines, Riquelme-Rodríguez et al. (2014) accounted for temporal variations in road humidity and developed an ALNS algorithm to solve the multi-objective arc routing problem aimed at minimizing sprinkling costs. Although some of these studies have considered aspects such as multi-depots, mixed graphs, and temporal changes, none of them have yet considered the response to real-time new demands during sprinkling operations. Thus, this study on routes planning of sprinklers on urban streets will address this issue.

The Dynamix-CARP proposed in this paper is a variant of the CARP, which has wide applications in many different fields (Mourão and Pinto, 2017). Bautista et al. (2008), Mourão and Almeida (2000), and Martins et al. (2015) established CARP models to describe the waste collection process in different cities, respectively. For waste collection, Willemse and Joubert (2016) introduced the mixed capacitated arc routing problem under time restrictions with intermediate facilities (MCARPTIF) and proposed three mechanisms to accelerate the solving process of large-scale problems (Willemse and Joubert, 2019). For the practical application of urban street services, Lai et al. (2020) conducted an investigation and proposed a tabu search algorithm for the split delivery capacitated arc routing problem with time windows (SDCARPTW), considering the properties of the split-delivery structure. To address urban traffic pollution more effectively, Cao et al. (2021) formulated a multi-objective CARP with four objective functions, namely total cost, makespan, carbon emission, and load utilization rate and developed a memetic algorithm (MA) to solve it. Khajepour et al. (2020) modeled agricultural harvest operations as a CARP and devised an ALNS algorithm for large-scale problems. These studies provide good examples on how to model practical problems as variants of the CARP.

Solution algorithms for the CARP and its variants have been extensively investigated and can be categorized into exact approaches and heuristic approaches. Due to the inherent NP-Hard properties of CARP, exact approaches are primarily applied to solve small-scale instances (Belenguer et al., 2015; Chen et al., 2016; Hintsch et al., 2021; Krushinsky and Van Woensel, 2015; Rivera et al., 2016). For large-scale CARP instances, heuristic and meta-heuristic methods have received significant attention (Prins, 2015). Zhang et al. (2017) proposed a MA that incorporated a new routing decomposition (RD) operator for the periodic capacitated arc routing problem (PCARP) and empirically verified the efficacy of MARD in solving large-scale PCARP instances. For the CARP, Arakaki and Usberti (2019) introduced an efficiency rule that selects the most promising edges to traverse next and proposed a new path-scanning heuristic, which outperformed all previous path-scanning heuristics. For the multi-trip capacitated arc routing problem (MCARP), Tirkolaee et al. (2019) designed a highly efficient ant colony optimization (ACO) algorithm based on an improved max-min ant system (MMAS). Keenan et al. (2021) examined the time capacitated arc routing problem (TCARP) and introduced a strategic oscillation simheuristic algorithm as a solution strategy.

Among various heuristic and meta-heuristic algorithms, one prominent approach used for solving large-scale CARP is the ALNS algorithm, which has exhibited exceptional performance in tackling large-scale problems (Mara et al., 2022). ALNS algorithm was initially proposed by Ropke and Pisinger (2006) to solve the pickup and delivery problem with time windows. This algorithm comprises several operators that can either destroy or repair the solution. Each operator is assigned a weight to control how frequently a particular operator is employed during the search. These weights are dynamically adjusted in the search process, enabling the algorithm to adapt to the specific instance and search status (Gendreau and Potvin, 2010). Since then, ALNS-based algorithms have been widely applied in various domains (Chen et al., 2021; Khajepour et al., 2020; Ma et al., 2023; Zhang et al., 2021), and their effectiveness has been extensively demonstrated (Franceschetti et al., 2017; Liu et al., 2023; Sacramento et al., 2019).

Furthermore, researchers have explored various improvements to the ALNS algorithm to enhance its effectiveness and solution quality. Alinaghian and Shokouhi (2018) proposed a hybrid algorithm for the multi-depot multi-compartment vehicle routing problem (MDMCVRP) by incorporating the variable neighborhood search (VNS) into ALNS algorithm and showed that the proposed hybrid algorithm outperformed the original algorithm through experiments. Mofid-Nakhaee and Barzinpour (2019) devised a hybrid ALNS algorithm integrated with the whale optimization algorithm (WOA) to solve the multi-compartment capacitated arc routing problems with intermediate facilities (MCCARPIF). Pitakaso et al. (2019) developed four new formulas to calculate the probability of accepting worse solutions in ALNS algorithm and found that the parabolic function achieved superior results. Cai et al. (2022) developed a hybrid algorithm combining ALNS algorithm with tabu search (TS) for the electric vehicle relocation problem and verified the competitiveness of their algorithm in terms of solution quality and computational efficiency. Moreover, the combination of TS and ALNS algorithm has emerged as a common approach for improving the algorithm's performance (He et al., 2020; Kır et al., 2017; Žulj et al., 2018).

In conclusion, the literatures on CARP and its variants are quite extensive. However, none of them can be directly used to solve the Dynamix-CARP and very few considered real-time responses to new demands during the sprinkling operations. In this paper, we proposed the Dynamix-CARP model for routes planning of sprinklers on urban streets and designed an improved ALNS algorithm to solve the problem. This study can serve as an important reference for operation companies to optimize the routes of sprinklers.

## 3. Problem formulation

### 3.1 Problem description

The Dynamix-CARP investigated in this paper arises from the sprinkling operations on urban streets. In sprinkling operations, sprinklers aim to mitigate FRD's adverse effects and improve air quality by sprinkling water on dusty streets within the study area.

As the typical water refilling time exceeds 30 minutes, the service time that a sprinkler with a full tank of water can provide is regarded as a planning cycle. For each cycle, streets that require sprinkling are identified as fixed demands and each street called demand edge. When the street has separating barriers in the middle, the sprinklers need to operate in both directions and the demand edge is regarded as directional. When the street has not have separating barriers, the sprinkler can operate in either direction and the demand edge is non-directional.

In the study area, there may be multiple depots. When the sprinkling operations start, the water-filled sprinklers depart from the depots simultaneously and follow the predetermined routes to sprinkle water on designated streets. The required number of sprinklers is roughly estimated by dividing the total length of the streets that need water-sprinkling by the length of the streets that one sprinkler can cover with a full tank of water. The working time of each sprinkler should not vary significantly so that all sprinklers have roughly the same workload. Finally, upon completion of the operations, these sprinklers will return to their departure depots to refill water.

During sprinkling operations, there may be real-time new sprinkling demands (real-time FRD pollution instance). In such cases, sprinklers need to serve the corresponding streets within a certain time window, by taking into consideration numerous factors such as tank load, travel distance, and deviation from designated route.

## 3.2 Notation

The Dynamix-CARP proposed in this paper is formulated based on a mixed graph $G = (V, L)$, as depicted in **Figure 1**. $V$ is a set of vertices that includes the depot set $V_o = \{O_1, O_2, \ldots, O_o\}$ and the node set $V_p = \{P_1, P_2, \ldots, P_p\}$, i.e., $V = V_o \cup V_p$. $L$ is the set of links in graph G including the directional arc set $A = \{<i, j> | i, j \in V\}$ and the non-directional edge set $E = \{(i, j) | i, j \in V\}$, i.e., $L = A \cup E$. For notational convenience, the ($n$th) link from node $i$ to node $j$ is denoted by $[i, j]$, $[i, j] \in$ L (or $l_n, l_n \in L, l_n = [i, j]$), so the directional arc set and the non-directional edge set can be represented as $A = \{[i, j] | i, j \in V\}$ and $E = \{[i, j] \cup [j, i] | i, j \in V\}$, respectively. Each link $[i, j]$ (or $l_n$) is associated with a distance $D_{ij}$ (or $D_n$), where $D_{ij} = D_{ji}$. Based on the distances of the links and their connection relationships, we utilize the Floyd-Warshall shortest-path algorithm to compute the shortest paths among all links (Willemse, 2016) and obtain the shortest travel distances between any two links, i.e., $D'_{mn} = D'_{nm}$, $l_n, l_m \in L$. Furthermore, the depots are modeled by incorporating them in $L$ as a set of virtual links with zero distance, i.e., $L_{V_o} = \{[i, i] | i \in V_o\}, L_{V_o} \in L$ and $D_n = 0, l_n \in L_{V_o}$.

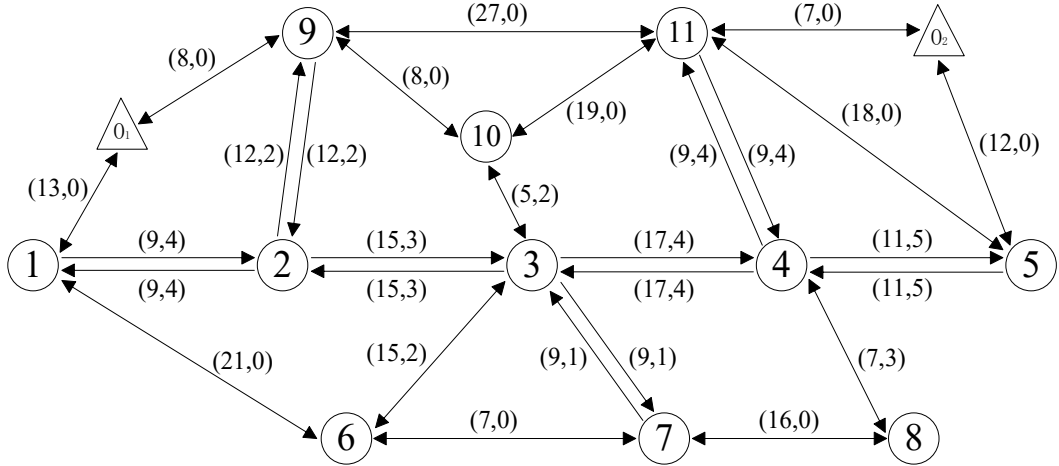

**Figure 1. An example of the mixed graph.**

The set of links $L$ can be further classified into two categories: the demand link set $L_R$ and the non-demand link set $L_N$. The demand link $[i, j] \in L_R$ has a positive demand ($d_{ij} > 0$) and must be served. The non-demand link $[i, j] \in L_N$ has a zero demand ($d_{ij} = 0$) and does not require sprinkling. Although these non-demand links do not need to be served, they are also part of the graph and contribute to the graph's connectivity. Furthermore, each link is allowed to be traversed multiple times. The corresponding $A_R$, $E_R$, $A_N$, and $E_N$ denote the demand arc set, demand edge set, non-demand arc set, and non-demand edge set, respectively. Demand arcs need to be sprinkled along the particular direction and demand edges can be sprinkled along either direction.

Assuming the fleet of sprinklers consists of an adequate number of identical sprinklers with the water-tank of capacity $Q$. The required number of sprinklers $m$ is determined by the total length of the streets that need water-sprinkling and the set of sprinklers is $K$, where $K = \{k_1, k_2, ..., k_m\}$. The $m$ sprinklers depart from the depots simultaneously at the departure time $t_0 = 0$. $T_n^k$ is the moment when sprinkler $k \in K$ starts to service the demand link $l_n \in L_R$ and $T_n^k = t_0, l_n \in L_{V_o}$. When receiving real-time sprinkling demands, the set $L_R'$ represents the demand links that have not been sprinkled, where $L_R' = E_R' \cup A_R'$ and the set of new demand links that need to be added to the routes is $L_R^{new}$, where $L_R^{new} = E_R^{new} \cup A_R^{new}$ and $L_R^{new} \in L$. The time window for the starting time of sprinkling on the new demand link $l_n \in L_R^{new}$ is $[a_n, b_n]$. Violation of the time window will result in a penalty and the coefficient related to the penalty is $\delta$.

The travel speeds when the vehicle is sprinkling and not sprinkling are $v_R$ and $v_N$, respectively. Usually, $v_R$ is lower than $v_N$. The allowed maximum difference between working time of sprinklers is $T_{dif}$. Additionally, each sprinkler $k \in K$ corresponds to a specific route $s_k$ within the solution $S$, i.e., $s_k \in S$, where $S = \{s_1, s_2, ..., s_m\}$.

For each route $s_k \in S$, let:

- $X_{ij}^k$ be 1 if sprinkler $k \in K$ serves the demand link $[i, j] \in L_R$ from $i$ to $j$, and 0 otherwise;
- $Y_{mn}^k$ be 1 if sprinkler $k \in K$ departs from the corresponding depot and firstly serves the demand link $l_n \in L_R$ ( $l_m \in L_{V_o}$ ), or proceeds to service demand link $l_n \in L_R$ after completing service demand link $l_m \in L_R$, and 0 otherwise.

### 3.3 Two-stage model

The response to new sprinkling demands in a timely manner makes the sprinkling routes planning a dynamic routing problem. In this study, we propose a two-stage model to generate the planned routes and the adjusted routes separately.

The first stage of the model is to generate the fixed routing plan based on the initial fixed sprinkling demands. The objective in this stage is to minimize the total travel distance.

Stage-1:

$$Min \sum_{k \in K} \sum_{[i,j] \in L_R} X_{ij}^k D_{ij} + \sum_{k \in K} \sum_{l_m, l_n \in L_R \cup L_{V_o} : l_n \neq l_m} Y_{mn}^k D_{mn}' . \tag{1}$$

s.t.:

$$X_{ij}^k = 0, \forall k \in K, \forall [i, j] \in L_N, \tag{2}$$

$$X_{ij}^k = \sum_{l_n \in L_R \cup L_{V_o} : l_n \neq l_m} Y_{mn}^k, \forall k \in K, \forall l_m = [i, j] \in L_R, \tag{3}$$

$$\sum_{l_m \in L_R \cup L_{V_o} : l_m \neq l_n} Y_{mn}^k = \sum_{l_m \in L_R \cup L_{V_o} : l_m \neq l_n} Y_{nm}^k, \forall k \in K, \forall l_n \in L_R, \tag{4}$$

$$\sum_{k \in K} X_{ij}^k + \sum_{k \in K} X_{ji}^k = 1, \forall [i,j] \in E_R, \tag{5}$$

$$\sum_{k \in K} X_{ij}^k = 1, \forall [i,j] \in A_R, \tag{6}$$

$$(m-1)Q < \sum_{k \in K} \sum_{[i,j] \in L_R} X_{ij}^k d_{ij} \leq mQ, \tag{7}$$

$$0 < \sum_{[i,j] \in L_R} X_{ij}^k d_{ij} \leq Q, \forall k \in K, \tag{8}$$

$$\sum_{l_m \in L_{V_o}} \sum_{l_n \in L_R} Y_{mn}^k = 1, \forall k \in K, \tag{9}$$

$$\sum_{l_n \in L_R} Y_{mn}^k = \sum_{l_n \in L_R} Y_{nm}^k, \forall k \in K, \forall l_m \in L_{V_o}, \tag{10}$$

$$T_n^k \geq T_m^k + \frac{D_m}{v_R} + \frac{D'_{mn}}{v_N} - M(1 - Y_{mn}^k), \forall k \in K, \forall l_m \in L_R \cup L_{V_o}, \forall l_n \in L_R: l_n \neq l_m, \tag{11}$$

$$T_{work}^k = \sum_{[i,j] \in L_R} \frac{X_{ij}^k D_{ij}}{v_R} + \sum_{l_m, l_n \in L_R \cup L_{V_o}: l_m \neq l_n} \frac{Y_{mn}^k D'_{mn}}{v_N}, \forall k \in K, \tag{12}$$

$$\max_{k \in K}(T_{work}^k) - \min_{k \in K}(T_{work}^k) \leq T_{dif}, \tag{13}$$

$$\mu_m^k - \mu_n^k + M(Y_{mn}^k - 1) + 1 \leq 0, \forall k \in K, \forall l_m, l_n \in L_R \cup L_{V_o}, l_m \neq l_n, \tag{14}$$

$$X_{ij}^k \in \{0,1\}, \forall k \in K, \forall [i,j] \in L, \tag{15}$$

$$Y_{mn}^k \in \{0,1\}, \forall k \in K, \forall l_m, l_n \in L_R \cup L_{V_o}, l_m \neq l_n, \tag{16}$$

$$T_n^k \geq 0, \forall\, k \in K, \forall l_n \in L_R \cup L_{V_o}, \tag{17}$$

$$\mu_n^k \geq 0, \forall\, k \in K, \forall l_n \in L_R \cup L_{V_o}. \tag{18}$$

The objective function, represented by equation **(1)**, is to minimize the total travel distance, composed of the distance traveled for both sprinkling and not sprinkling. Constraint **(2)** indicates that non-demand links do not require sprinkling. Constraint **(3)** and Constraint **(4)** make sure the continuity of each sprinkler's route. Constraint **(5)** ensures that each demand edge can be sprinkled along either direction but only once. Constraint **(6)** guarantees each demand arc is covered by one sprinkler along the specified direction. Constraint **(7)** determines the number of sprinklers. Constraint **(8)** ensures that the total water consumption for each route does not exceed the capacity of the sprinkler's water-tank. Constraint **(9)** means that each sprinkler must start from a depot and cannot return to the depots in the midway. Constraint **(10)** guarantees that each sprinkler eventually returns to the departure depot. Constraint **(11)** updates the start time of sprinkling for each demand link. Constraint **(12)** and Constraint **(13)** make sure the difference in working time among the sprinklers is within a certain limit. Constraint **(14)** aims to eliminate isolated subtours. Constraints **(15)-(18)** specify the ranges of values for corresponding variables.

In the second stage, some demand links have already been served, necessitating route adjustments to promptly respond to the real-time new sprinkling demands. When making adjustments, the sprinklers' current positions, the amount of remaining water, and the constraints specified before need to be considered. Notably, the demand links for this stage include both new sprinkling demands that need to be added to the routes and the initial fixed demands that have not been served.

Stage-2:

$$Min \sum_{k \in K} \sum_{[i,j] \in L_R^{new} \cup L_R} X_{ij}^k D_{ij}$$

$$+ \sum_{k \in K} \sum_{l_m, l_n \in L_R^{new} \cup L_R \cup L_{V_o} : l_n \neq l_m} Y_{mn}^k D_{mn}^{'} + \delta \sum_{k \in K} \sum_{l_n \in L_R^{new}} \max^2(0, T_n^k - b_n). \qquad (19)$$

s.t. **(2)-(18)**, and

$$T_n^k \geq a_n, \forall\, k \in K, \forall l_n \in L_R^{new}. \qquad (20)$$

As the new sprinkling demands are usually urgent, there are time window restrictions for such operations. The objective function expressed in Equation **(19)** minimizes the sum of the total travel distance and the penalty cost for not serving the new demands in time (beyond the time window). Constraint **(20)** ensures that the start time of sprinkling for the demand link with time window is not earlier than the start time of the corresponding time window.

## 4. Solution algorithm

In this section, we present a customized algorithm for solving the Dynamix-CARP based on the improved ALNS. **Algorithm 1** provides the pseudo-code for the proposed algorithm, which starts by constructing an initial solution $S$ (described in **Algorithm 2**) containing all the fixed demand links. Subsequently, the improved ALNS algorithm (described in **Algorithm 3**) is applied to obtain an optimized solution $S^*$. Then the sprinklers travel along the routes specified in $S^*$.

**Algorithm 1** Pseudo-code of the algorithm

1: $S \leftarrow InitSol\ (L_R, V_o, Q, T_{dif})$ // Algorithm 2
2: $S^* \leftarrow ALNS(S)$ // Algorithm 3
3: Return $S^*$
4: **while** receiving real-time sprinkling demands **do**
5:     Construct the set of demand links that need to be added to the current routes: $L_R^{new}$
6:     $S^{'} \leftarrow Update(S^*)$ // Remove sprinkled links
7:     **for** $i$ in $L_R^{new}$ **do**
8:         Add time window to $i$
9:         **if** $i$ in $S^{'}$ **then**
10:             Remove $i$ from $S^{'}$
11:         **end if**
12:     **end for**
13:     $S_D \leftarrow Repair(S^{'}, L_R^{new})$ // Insert $L_R^{new}$ into $S^{'}$
14:     $S_D^* \leftarrow ALNS(S_D)$
15:     Return $S_D^*$
16: **end while**

During the operation, when receiving real-time sprinkling demands (line 4), the algorithm will determine the demand links that need to be added to the current routes (line 5) and remove the links that have already been served in the solution $S^*$, resulting in a new solution $S'$ (line 6). The next step is to set the time window for the corresponding demand links (line 8) and if these demand links are already in $S'$, they will be removed from $S'$ (line 10). Subsequently, the demand links that need to be added to the current routes are incorporated into the $S'$ using the repair operator, yielding the initial solution $S_D$ for subsequent sprinkling operations (line 13). Finally, $S_D$ is further optimized by the improved ALNS algorithm to generate the final solution $S_D^*$ (line 14).

The following subsections provide a detailed description of each part of the algorithm.

### 4.1 Initial solution

The initial solution is generated using the greedy insertion method described in **Algorithm 2**.

**Algorithm 2** *InitSol*
─────────────────────────────────────────────
1: **Input:** $L_R$, $V_o$, $Q$, $T_{dif}$
2: $R_s \leftarrow L_R$ // All demand links are initially unserved
3: $m \leftarrow int\left(\frac{total\ demand}{Q} + 0.5\right)$ // Sprinklers number
4: $L_1 \leftarrow random\,(L_R)$, $s_1 \leftarrow [L_1]$, remove $L_1$ from $R_s$, $A_s \leftarrow [L_1]$, $r \leftarrow 1$, $now_1 \leftarrow j : L_1 : [i, j]$
5: **while** $r < m$ **do** // While there are empty routes
6:     $++r$
7:     $L_r \leftarrow \arg\max\limits_{L}\left\{\sum\limits_{A \in A_s} distance\,(L, A)\,, L \in R_s\right\}$
8:     $s_r \leftarrow [L_r]$, remove $L_r$ from $R_s$, $A_s \leftarrow A_s \cup [L_r]$, $now_r \leftarrow j : L_r : [i, j]$
9: **end while**
10: **while** $R_s$ is not $\emptyset$ **do** // While there are demand links unassigned
11:     $f \leftarrow 1$
12:     **while** $f \leq m$ **do** // Assign links to the routes in turn
13:         $L^* \leftarrow \arg\min\limits_{L}\left\{distance\,(L, now_f)\,, L \in R_s\right\}$
14:         **if** $\left(\sum demand\,(s_f) + demand\,(L^*)\right) > Q$ **then**
15:             $++f$
16:         **else**
17:             Add $L^*$ to the end of the $s_f$, remove $L^*$ from $R_s$, $now_f \leftarrow j : L^* : [i, j]$, $++f$
18:         **end if**
19:     **end while**
20: **end while**
21: $S \leftarrow [s_1, \ldots, s_m]$
22: **for** $s_i$ in $S$ **do**
23:     $p_i \leftarrow \arg\min\limits_{p}\left\{\sum\limits_{L \in s_i} distance\,(L, p)\,, p \in V_o\right\}$
24:     Insert $p_i$ into $s_i$ at the position that minimizes the total travel distance of $s_i$
25: **end for**
26: **while** $\max\{work\ time\,(s_f)\,, s_f \in S\} - \min\{work\ time\,(s_f)\,, s_f \in S\} > T_{dif}$ **do**
27:     $r_1 \leftarrow \arg\max\limits_{s_f}\{work\ time\,(s_f)\,, s_f \in S\}$
28:     $r_2 \leftarrow \arg\min\limits_{s_f}\{work\ time\,(s_f)\,, s_f \in S\}$
29:     $L^{'} \leftarrow \arg\max\limits_{L}\{distance\,(p_{r_1}, L)\,, L \in r_1\}$
30:     Move $L^{'}$ from $r_1$ to $r_2$ at the position that minimizes the total travel distance of $r_2$
31: **end while**
32: Return $S$
─────────────────────────────────────────────

The initial solution begins with $m$ empty routes, where $m$ is the number of sprinklers engaged in this operation. A demand link $L_1$ is randomly selected as the link that requires service by the first sprinkler (route). Subsequently, the $m - 1$ demand links that are farthest from $L_1$ are selected as the links that need to be sprinkled by the remaining sprinklers. The endpoint of each selected demand link becomes the current position of the corresponding sprinkler. In the next step, the demand link that is the closest to the sprinkler's current position is added to the route of each sprinkler and the current position of the sprinkler is updated accordingly. If the total water consumption of a sprinkler reaches the capacity of the water-tank $Q$, this sprinkler is skipped while the remaining demand links are assigned to the other sprinklers until all the demand links have been assigned to the sprinklers.

Since there are multiple depots in the study area, it is necessary to determine the appropriate depot for dispatching sprinklers. In this approach, the depot closest to the route is selected and added to the route at appropriate location to minimize the corresponding sprinkler's travel distance.

Finally, it is essential to limit the difference in working time among the sprinklers to ensure that all sprinklers have roughly the same workload. If the difference between the longest and the shortest working time exceeds $T_{dif}$, the demand link that is farthest from the depot associated with the route of the longest working time is moved to the route of the shortest working time. This adjustment is iteratively executed until the difference in working time of the sprinklers does not exceed $T_{dif}$.

## 4.2 Adaptive large neighborhood search

The improved ALNS algorithm builds upon the work of Ropke and Pisinger (2006), with several modifications to suit the specific characteristics of the Dynamix-CARP. Firstly, specific destroy and repair operators are designed to modify the allocation and sequence of the demand links. Secondly, an acceleration mechanism is introduced to reduce the search domains during the repair process with large computation. Thirdly, the algorithm employs a tabu-list to record the search process and to avoid duplicate searches. Fourthly, perturbations are strategically added to expand the search space and worse solutions are considered when the search process fails to find improved solutions after multiple iterations. These modifications collectively enhance the efficiency and effectiveness of the improved ALNS algorithm in finding a better solution.

The general structure of the improved ALNS algorithm is presented in **Algorithm 3**, with its components described in detail later. In this algorithm, $\Omega^-$ and $\Omega^+$ represent the sets of destroy and repair operators, respectively. $\rho$ and $scores$ record the weights and scores of each destroy-repair pair, respectively. The first step is the variable initialization, followed by the optimization process performed in an iterative way. In each iteration, a destroy-repair pair is selected using the roulette wheel mechanism based on their weights and the tabu-list. The selected pair is applied to the current solution $S_{curr}$ to generate a new solution $S_{new}$. The relevant parameters are updated based on the values of $S_{new}$, $S_{curr}$, and the optimal solution $S_{best}$. When $F(S_{new}) > F(S_{curr})$, the acceptance criterion described by Ropke and Pisinger (2006) for simulated annealing is applied to determine whether $S_{new}$ could be accepted. In every $\psi$ iterations, the weight of each destroy-repair pair is updated according to the scores.

---
**Algorithm 3** *ALNS*
---
1: **Input:** A feasible solution $S$, $\Omega^-$, $\Omega^+$, $T_0$
2: $S_{best} \leftarrow S_{curr} \leftarrow S$, $Iter \leftarrow 1$, $nonIter \leftarrow 1$, $T \leftarrow T_0$
3: **Initialize:** $\rho \leftarrow (1, ..., 1)$, $scores \leftarrow (0, ..., 0)$, $tabulist \leftarrow \phi$
4: **repeat**
5:     Select a destroy-repair pair $(d, r)$ using the roulette wheel mechanism, $d \in \Omega^-$, $r \in \Omega^+$
6:     $S_{new} \leftarrow r(d(S_{curr}, tabulist))$
7:     Update $T$, $tabulist$
8:     $++ Iter$
9:     **if** $F(S_{new}) < F(S_{best})$ **then**
10:         $S_{best} \leftarrow S_{curr} \leftarrow S_{new}$, $nonIter \leftarrow 1$
11:     **else if** $F(S_{new}) < F(S_{curr})$ or $accept(S_{new}, S_{curr}, T)$ **then**
12:         $S_{curr} \leftarrow S_{new}$, $++ nonIter$
13:     **else**
14:         $++ nonIter$
15:     **end if**
16:     Update the score of the operator pair
17:     **if** $nonIter > \phi_1$ **then**
18:         **if** $nonIter$ mod $\phi_2 == 0$ **then**
19:             $S'_{curr} \leftarrow Shake(S_{curr}, n)$
20:             **if** $F(S'_{curr}) < F(S_{best})$ **then**
21:                 $S_{best} \leftarrow S_{curr} \leftarrow S'_{curr}$, $nonIter \leftarrow 1$
22:             **else if** $F(S'_{curr}) < F(S_{curr})$ or $accept(S'_{curr}, S_{curr}, T)$ **then**
23:                 $S_{curr} \leftarrow S'_{curr}$
24:             **else**
25:                 $S_{curr} \leftarrow random(S_{curr}, S_{before})$
26:             **end if**
27:         **else**
28:             $S_{curr} \leftarrow random(S_{curr}, S_{before})$
29:         **end if**
30:     **end if**
31:     **if** $Iter$ mod $\psi == 0$ **then**
32:         Update the $\rho$ and reset the $scores \leftarrow (0, ..., 0)$
33:     **end if**
34: **until** $nonIter > M_1$ or $Iter > M_2$
35: Return $S_{best}$
---

If the optimal solution is not updated after $\Phi_1$ iterations, the current solution or a worse solution is randomly selected as the current solution for the next iteration. Moreover, a perturbation is performed on the current solution once every $\Phi_2$ iterations. If a better solution or an acceptable solution is still not obtained after the perturbation, the current solution or the worse solution is randomly selected as the next current solution.

Finally, if the optimal solution remains unchanged for $M_1$ iterations or when the maximum number of iterations $M_2$ is reached, the algorithm terminates.

### 4.2.1 Adaptive weighting mechanism

In the proposed ALNS algorithm, each destroy-repair pair is assigned a weight $\rho_i$, and these weights are dynamically adjusted based on the performance of each pair. Initially, all destroy-repair pairs have the same weight (i.e., $\rho_{i1} = 1$, $i = 1,2,\dots,p$), which ensures that each pair has an equal probability of being selected, with a selection probability of:

$$P_{i1} = \frac{\rho_{i1}}{\sum_{i=1}^{p} \rho_{i1}}. \tag{21}$$

The search process is divided into phases of 100 iterations each. At the end of each phase, the weight of each destroy-repair pair is adjusted based on the performance of the pair in that phase:

$$\rho_{in} = \tau\rho_{in-1} + \frac{score_i}{time_i}, \tag{22}$$

where $\tau \in [0,1]$ controls the effect of the previous weights on the current weights. The $score_i$ and $time_i$ are the scores of destroy-repair pair $i$ and its usage times in the current stage, respectively.

Following the suggestion of Ropke and Pisinger (2006), three scores ($sco_1, sco_2, sco_3$) are defined to characterize the performance of each destroy-repair pair when the newly generated solution is accepted. Here, $sco_1 < sco_2 < sco_3$, and they are used for destroy-repair pair $i$ in cases where $F(S_{curr}) < F(S_{new})$, $F(S_{best}) < F(S_{new}) < F(S_{curr})$, and $F(S_{new}) < F(S_{best})$, respectively.

To prevent premature convergence, the lower limit and the upper limit of each destroy-repair pair's weight are set to $a_{min}$ and $a_{max}$, respectively.

*4.2.2 Destroy operation*

After the solution is obtained, since adjusting the position of some demand links within the solution may improve the current solution, the destroy operation is applied to remove some demand links from the current solution that could contribute to a larger objective function. The proposed ALNS algorithm utilizes five destroy operators.

(1) Random removal

This operator randomly selects $\Gamma$ demand links and removes them from their current positions in the routes. While this kind of removal might lead to a worse solution, it could help diversify the search and potentially prevent the solution from being a local optimum.

(2) Worst removal

The purpose of this operator is to identify and remove $\Gamma$ demand links that could increase the objective function significantly. Therefore, if removing a demand link from the current solution results in an evident decrease in the objective function, the corresponding link will be removed to improve the solution.

(3) Non-adjacent links removal

This operator calculates the distances between each demand link and its upstream and downstream demand links in the solution, respectively. The sum of both distances is the non-adjacent distance. The longer the non-adjacent distance of a demand link is, the less desirable it is to include the corresponding demand link sequence in the solution. Thus, this operator identifies and removes the $\Gamma$ demand links with the longest non-adjacent distances from the current solution.

(4) Farthest depot removal

Since the sprinklers must return to the departure depot after completing their operation, the ideal solution should let sprinklers serve the demand links that are close to the departure depot. If the demand links served by a sprinkler are far away from the sprinkler's departure depot in the current solution, it may be beneficial to assign those demand links to sprinklers departing from other depots.

(5) Time-related removal

This operator removes $\Gamma$ demand links from the current solution that cause a penalty cost for sprinkling beyond the time window. If the number of this kind of demand links is less than $\Gamma$, all such links are removed and the remaining demand links are randomly removed until the number of removed links equals to $\Gamma$. It is important to note that the time window restrictions are only effective when adjusting the sprinkling routes. During the fixed routes planning, this operator is equivalent to the random removal operator.

### 4.2.3 Repair operation

The repair operation is responsible for generating a new feasible solution by inserting the removed demand links into the destroyed solution at appropriate positions. The proposed ALNS algorithm utilizes five repair operators.

(1) Greedy insertion

The greedy insertion is a commonly used exact repair method in routing problems. It aims to minimize the insertion cost by identifying the positions where inserting the demand links causes the least increase in the objective function.

(2) Noise greedy insertion

Noise greedy insertion extends the greedy insertion by introducing a noise term, which is designed to decrease with the number of iterations. The formula for the noise term is as follows:

$$\Delta n = \mu \varphi z d_{max}, \tag{23}$$

where $d_{max}$ is the maximum non-adjacent distance between demand links, $\mu$ is a parameter controlling the magnitude of the noise, $\varphi = 1 - \frac{iter}{all\_iter}$ regulates the fluctuation range of the noise decreasing as the number of iterations progresses, and $z \in [0,1]$ is a random number.

(3) Regret insertion

The greedy insertion may be quite myopic, making it challenging for later inserted demand links to find favorable inserting positions due to the previous insertion of other demand links, which reduces some of the available inserting positions. To address this, regret insertion is introduced as a look-ahead method and used in conjunction with greedy insertion.

The 2-regret criterion is employed in this process. It is implemented by first calculating the regret value, which is the difference of the objective function between the optimal position and the sub-optimal position to insert each demand link, and then preferentially inserting the demand link with the maximum regret value to the optimal position.

(4) Noise regret insertion

Similar to greedy insertion, regret insertion is also an exact repair method. Therefore, the noise term is also incorporated into the regret insertion to calculate the insertion cost, aiming to increase the diversity of the solution.

(5) Time window insertion

Time window insertion is specifically designed for demand links with the sprinkling time window. Sprinkling beyond the time window can lead to a significant increase in the objective function, so the demand links that have the sprinkling time window will be inserted into positions that do not go beyond the time window while minimizing the increase in the objective function.

During the phase of fixed routes planning, after the repair operation, the depot of the original route may not be the one that minimizes the objective function. The algorithm needs to check whether the depot of each route after the repair operation corresponds to the optimal solution. If not, the depot of the route will be modified accordingly.

Additionally, this paper introduces an accelerated approach to reduce the unnecessary search domains for repair operators. When searching for the optimal inserting location of each demand link, if the link is very close to a specific depot $i$ and significantly farther away from other depots, then the ideal solution would be letting a sprinkler departing from depot $i$ serve this demand link. Therefore, there is no need to search for the routes of sprinklers departing from depots that are far away from this demand link. Thus, the search domain is significantly reduced and the efficiency of the algorithm is improved, which is particularly useful when the number of demand links is large.

Furthermore, it is essential to note that certain exact destroy and repair operators, such as worst removal, non-adjacent links removal, farthest depot removal, and time-related removal, along with greedy insertion and regret insertion, generate the same new solution each time when they are applied to the same solution. To avoid redundant searches, we employ a tabu-list to record the operations performed on the solutions using these combinations of destroy operators and repair operators.

*4.2.4 Perturbation*

To introduce diversification in the search process, perturbations are integrated into the algorithm. If the optimal solution is not updated after $\Phi_1$ iterations, the current solution undergoes perturbation every $\Phi_2$ generations. During the perturbation stage, a higher proportion $\Gamma_{shake}$ of destruction is applied to the current solution. Subsequently, two new solutions are generated using two repair methods, greedy insertion and regret insertion, respectively. The better of the two new solutions is chosen. When the chosen solution is not superior to the current solution, the principle of simulated annealing is used to decide whether it should be accepted or not. When the chosen solution is not better than the current solution and is not accepted, either the current solution or the worse solution is randomly chosen as the new current solution.

### 4.3 Dynamic sprinkling route adjustments

During sprinkling operations, it is necessary to adjust the sprinkling routes to respond to the real-time new demands caused by a combination of factors such as construction related truck activities, surge in traffic volume, and dry and windy weather. Typically, such new demands are related to a series of neighboring road sections.

Regarding the adjustment of sprinkling routes, the first step is to determine the set of the demand links that need to be added to the current routes, denoted as $L_R^{new}$. If the new demand links are in the current routing plan but have not yet been served, promptly responding to real-time demands is equivalent to bringing the corresponding sprinkling tasks forward. For new demand links that are in the current routing plan and have already been sprinkled, if they were sprinkled recently (such as in the past 30 minutes), there is no need to serve them again. But if they were sprinkled a long time ago or are not in the current routing plan, they need to be served as soon as possible. In summary, new demand links that require timely sprinkling include those originally in the current routing plan but not yet sprinkled, those sprinkled a long time ago, and those not in the current routing plan. If the remaining water of the engaged sprinklers can meet the demands of both the original unserved demand links and the new demand links requiring timely sprinkling, the latter are considered as the demand links that need to be added to the routes. Otherwise, based on the remaining water, the algorithm randomly selects the links from the new demand links that require timely sprinkling as the demand links that need to be added to the current routes.

Subsequently, the algorithm adds the time windows to the demand links that need to be added to the routes, requiring that the sprinkling of these demand links needs to be executed within a specific time frame. It is important to note that, upon completing the sprinkling operations, the sprinklers should return to their respective departure depots. Consequently, the current position of each sprinkler serves as the starting point for its respective route, while the original end depots of the routes remain unchanged, generating the current state solution. Then the algorithm inserts the demand links that need to be added to the routes at appropriate positions within the current state solution and generates a new dynamic initial solution. This can be accomplished using any of the repair operators described in section **4.2.3** Repair operation.

Finally, through the optimization of the dynamic initial solution using the improved ALNS algorithm detailed in section **4.2** Adaptive large neighborhood search, the sprinkling routes scheme for the next stage is obtained.

## 5.  Cases study

To validate the operational feasibility and effectiveness of the proposed algorithm for Dynamix-CARP, we perform the cases study based on data of Yulin Block in Chengdu, China.

The algorithms are programmed using Python 3.8 and executed on a machine with AMD Ryzen 9 5950X 16-Core Processor.

## 5.1 Data description and parameter setting

Yulin Block is located in the center of Chengdu and has an area of 5.3 km$^2$. The complete road network is obtained from the Open Street Map as shown in **Figure 2(a)**. To accommodate the routing algorithm, we extract the road sections that sprinklers can traverse smoothly and serve effectively. The simplified road network is shown in **Figure 2(b)**, which consists of 210 vertices, 178 directional arcs, and 235 non-directional edges. The sprinklers employed by local authority have a maximum sprinkling distance of 20 km (with full tanks), travels at a speed of $v_R = 10 km/h$ when sprinkling, and $v_N = 30 km/h$ when not in operation. The maximum difference between sprinklers' working time, $T_{dif}$, is set to be 15 minutes.

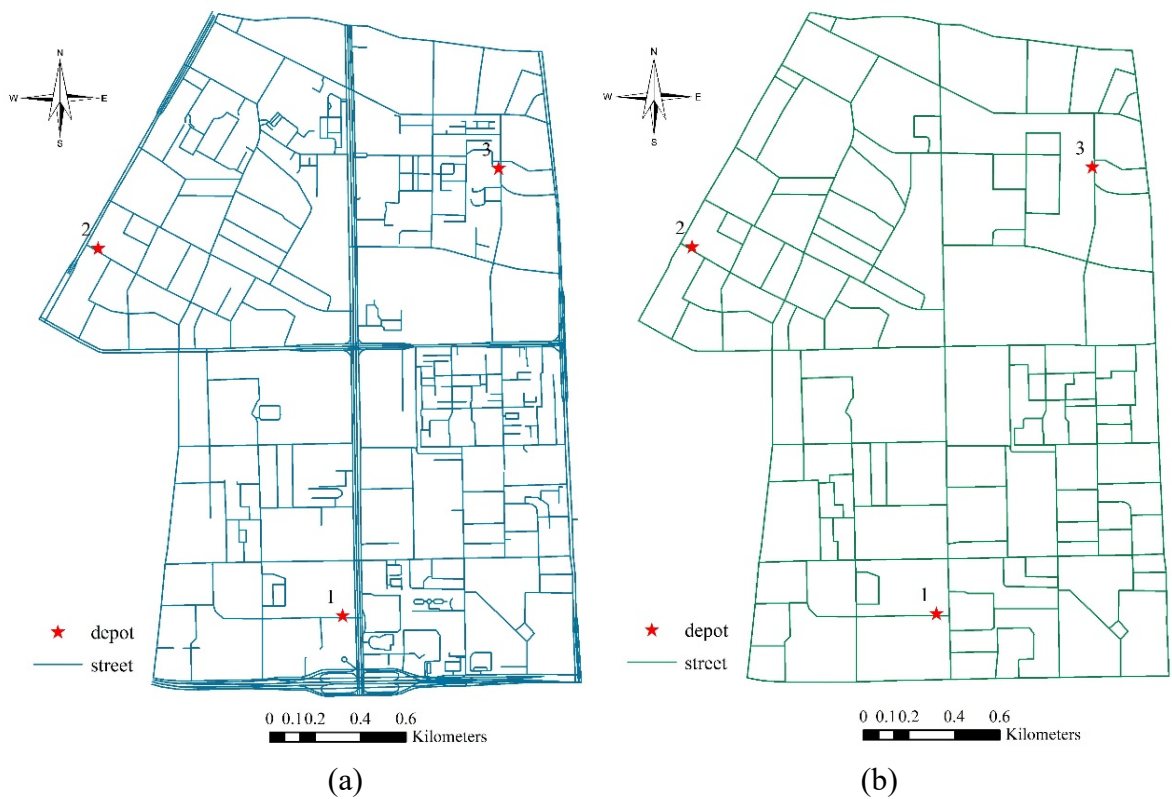

(a)                                    (b)

**Figure 2. (a) Full road network from Open Street Map; (b) simplified road network.**

For the convenience of testing algorithms, we specify that real-time new sprinkling demands were received 60 minutes after the start of the operations and these new demands should be sprinkled within 30 minutes upon receipt. The test dataset includes only real-time sprinkling demands once. For multiple real-time route adjustments, the same set of steps can be repeated accordingly.

The algorithm parameters are fine-tuned to balance solution quality and computation time based on preliminary experiments. The values used are compiled in **Figure 1**.

**Table 1 Parameter Values.**

| Parameter | Value | Parameter | Value |
|---|---|---|---|
| $\Gamma$ | 10% | $\psi$ | 100 |
| $\delta$ | 5 | $a_{min}$ | 1 |
| $\Phi_1$ | 200 | $a_{max}$ | 3 |
| $\Phi_2$ | 100 | $sco_1$ | 0 |
| $M_1$ | 1500 | $sco_2$ | 0.02 |
| $M_2$ | 3000 | $sco_3$ | 0.05 |

## 5.2 Algorithm performance

In practical scenarios, the FRD concentration varies significantly at different times, resulting in fluctuations in sprinkling demands. We propose three scenarios to assess the algorithm's performance under different levels of fixed sprinkling demands. For each scenario, we randomly selected specific road sections as new demands that require timely sprinkling and designed two scenarios for real-time route adjustments. Consequently, a total of six scenarios are formulated to evaluate the effectiveness of the proposed algorithm. The results are compared to those obtained using the ALNS algorithm. Both algorithms are run 10 times under each scenario to generate the optimal solution. The scenarios and their corresponding results are presented in **Table 2** and **Table 3**.

**Table 2 Solution without real-time route adjustments.**

| Scenarios | Algorithm | $A_R$ | $E_R$ | $Obj$ | $t(s)$ | $Obj(1min)$ | $Obj(5min)$ | $ave$ | $num$ |
|---|---|---|---|---|---|---|---|---|---|
| 1, 2 | Proposed algorithm | 142 | 3 | 30061 | 30.5 | - | - | 31348 | 2 |
| | ALNS | | | 30668 | 33.6 | - | - | 32020 | |
| 3, 4 | Proposed algorithm | 158 | 81 | 54001 | 316.2 | 55340 | 54001 | 54941 | 3 |
| | ALNS | | | 54023 | 695.5 | 56479 | 54111 | 55156 | |
| 5, 6 | Proposed algorithm | 178 | 172 | 76303 | 1607.8 | 78485 | 77420 | 76534 | 4 |
| | ALNS | | | 76282 | 2651.5 | 80358 | 78091 | 77933 | |

**Table 3 Solution with real-time route adjustments.**

| Scenarios | Algorithm | $A_R^{new}$ | $E_R^{new}$ | $Obj$ | $t(s)$ | $Obj(1min)$ | $Obj(5min)$ | $ave$ | $num$ |
|---|---|---|---|---|---|---|---|---|---|
| 1 | Proposed algorithm | 12 | 6 | 35662 | 9.5 | - | - | 36845 | 2 |
| | ALNS | | | 36374 | 14.2 | - | - | 37674 | |
| 2 | Proposed algorithm | 32 | 9 | 40188 | 17.9 | - | - | 40657 | |
| | ALNS | | | 41504 | 25.0 | - | - | 41687 | |
| 3 | Proposed algorithm | 26 | 5 | 56613 | 47.1 | - | - | 56636 | 3 |
| | ALNS | | | 56807 | 95.9 | 56807 | - | 57373 | |
| 4 | Proposed algorithm | 50 | 12 | 61072 | 99.7 | 61082 | - | 61735 | |
| | ALNS | | | 62156 | 219.1 | 63896 | - | 63715 | |
| 5 | Proposed algorithm | 40 | 10 | 83078 | 447.1 | 85487 | 83215 | 83918 | 4 |
| | ALNS | | | 83929 | 702.4 | 85883 | 84706 | 84438 | |
| 6 | Proposed algorithm | 54 | 22 | 89788 | 596.4 | 92045 | 90751 | 90069 | |
| | ALNS | | | 90854 | 844.6 | 95185 | 91722 | 91935 | |

In **Table 2 and Table 3**, $A_R$ and $E_R$ denote the number of fixed demand arcs and fixed demand edges, respectively. $A_R^{new}$ and $E_R^{new}$ indicate the number of new demand arcs and new demand edges that need to be added to the current routes, respectively. $Obj$, $Obj(1min)$, and $Obj(5min)$ stand for the values of the objective function for the final, 1-minute and, 5-minute solutions, respectively. $t$ represents the time taken to find the final solution and $ave$ represents the mean value of the objective function after running 10 times. Additionally, $num$ represents the number of sprinklers engaged.

The results of the experiments reveal a positive correlation between the computation time and the amount of demands. In instance 1 and instance 2 where the number of fixed demand links is low, both algorithms generate the fixed routing plans within 35 seconds and complete the route adjustments within 30 seconds. The difference in computation time between the two algorithms is marginal. For medium-sized instance 3 and instance 4, both algorithms generate the fixed routing plans within 15 minutes and complete the route adjustments within 5 minutes. The proposed algorithm demonstrates higher efficiency, reducing more than half of the computation time compared to the ALNS algorithm. In the case of larger-scale instance 5 and instance 6, the ALNS algorithm requires nearly 45 minutes to generate the fixed routing plans. Although the proposed algorithm significantly reduces the computation time, it still takes 25 minutes to obtain the fixed routing plans and approximately 10 minutes to adjust the fixed routes.

Despite the algorithms spend longer computation time for larger-scale instances, the proposed algorithm still manages to obtain relatively good solutions in a shorter time. An assessment was conducted on the algorithm's ability to find the relatively good solutions quickly, considering the solutions obtained when the algorithms ran for 1 minute and 5 minutes to generate the fixed routing plans and make route adjustments. For medium-sized instance 3 and instance 4, the objective value of the solution generated by the proposed algorithm within 1 minutes is just 2.48% larger than that of the final solution. For the case of route adjustments, the solutions generated by the proposed algorithm within 1 minute are similar to the final solutions, and outperform the final solutions obtained by the ALNS algorithm. For larger-scale instance 5 and instance 6, although the gaps between the shorter-time solutions and the final solutions are larger than that of instance 3 and instance 4, they remain acceptable. Overall, the proposed algorithm exhibits good efficiency and effectiveness in finding relatively good solutions in a reasonable amount of time.

Regarding the solution quality, the proposed algorithm demonstrates robustness after multiple runs, and the discrepancies between the mean values of the objective value for the 10 solutions and the objective value for the optimal solutions are minimal. In most instances, the proposed algorithm could generate superior solutions in a shorter time. In a few scenarios, the proposed algorithm may produce slightly inferior solutions due to its reduction of the search domain. However, the proposed algorithm's clear advantage lies in its short computation time. Overall, the acceleration mechanism and algorithm improvement strategies employed in this paper have been proved to be effective in enhancing the algorithm's performance.

## 5.3 Sensitivity analysis

To assess the impact of different parameters on the planning results, we conducted the sensitivity analyses using the parameters and data from scenario 4. For each case, we ran the algorithm ten times and selected the optimal solution.

### 5.3.1 Response time

At various stages of the sprinkling operation, the positions of the sprinklers and the demand links that have not yet been sprinkled are different. Consequently, the time when adjusting the routes (receiving the real-time demands) will influence the subsequent planning results. This impact is shown in **Figure 3**, where the response time is defined as the duration between the moment that the sprinklers depart from the depots and the moment that the real-time demands are received.

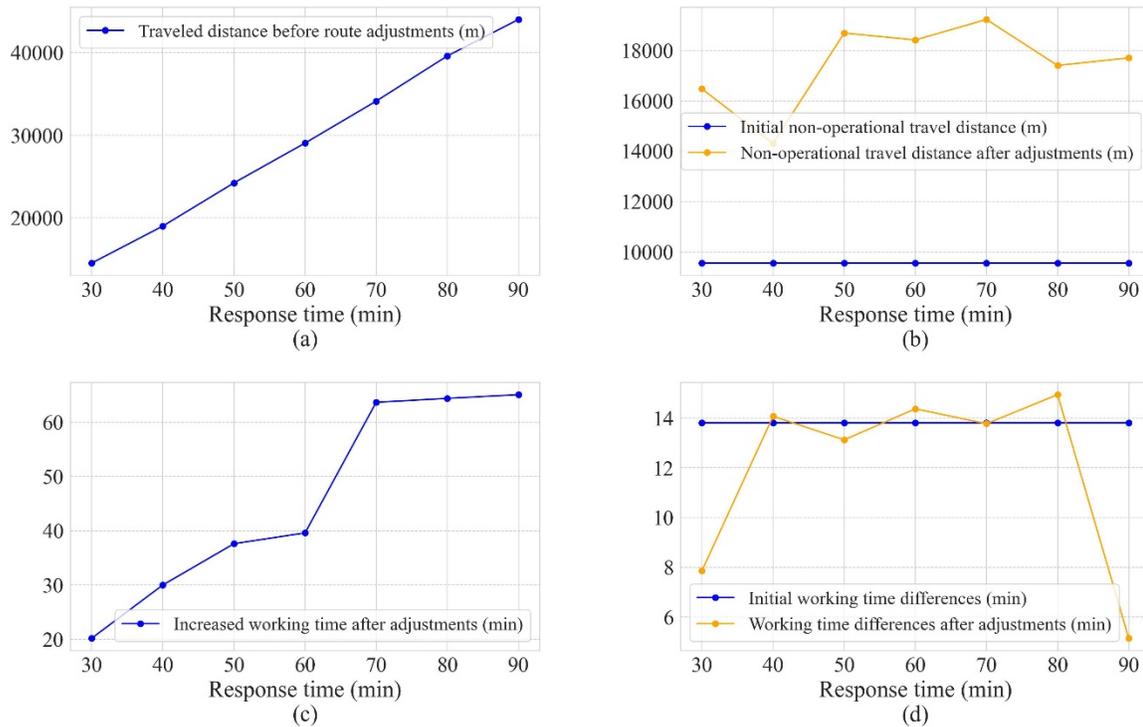

**Figure 3. Results under different response time**

**Figure 3**(a) illustrates that the sprinklers' traveled distance before route adjustments has a positive linear relationship with the response time, indicating that the workload distribution of the sprinklers remains relatively consistent over time.

**Figure 3**(b) reveals that route adjustments significantly increase the overall non-operational travel distance of the sprinklers. Once the response time surpasses 40 minutes, the sprinklers' non-operational travel distance after route adjustments reaches a higher level. This may be because when the route adjustments occur earlier, there are many unserved demand links. When sprinkling the new demand links, sprinklers can also sprinkle some adjacent unserved demand links, thus reducing the total non-operational travel distance. Conversely, when the response time exceeds 40 minutes, the number of unserved demand links diminishes, only a limited number of unserved demand links are close to the new demand links. Sprinklers prioritize serving the new demand links before attending to the remaining demand links, resulting in overlapping routes.

The working time of the sprinklers is determined by both the sprinkling distance and the non-operational travel distance. Due to the slower travel speed of the sprinklers during sprinkling compared to not sprinkling, the sprinkling distance has a more significant impact on the working time of the sprinklers than the non-operational travel distance. In **Figure 3**(c), the total working time of the sprinklers when real-time demands are received 60 minutes after the start of the operation has significantly increased compared to the previous moments. This increase may be attributed to some new demand links that require timely sprinkling were sprinkled a long time ago, thus necessitating repeat sprinkling. This results in a longer sprinkling distance and, consequently, an extended working time. Additionally, the later the route adjustments occur, the greater the total working time is.

The working time of the sprinklers differs, possibly due to variations in workload of each sprinkler and the distribution of demand links. **Figure 3**(d) depicts the variance in working time differences among the sprinklers concerning the response time, showing evident difference in the working time of the sprinklers within the fixed routing plans. When the real-time demands are received between 40 and 80 minutes after the start of the operation, the differences in working time among the sprinklers after route adjustments fluctuate around the difference in the fixed routing plans. This indicates the significant influence of the fixed routes on the differences in the working time of the sprinklers after route adjustments. Additionally, it is noteworthy that the route adjustments, whether occurring very early or very late, can effectively reduce the differences in working time among the sprinklers.

*5.4.2 Response time window*

This study requires that the new demand links should be sprinkled within a 30-minute time window when responding to the real-time new sprinkling demands. Delayed responses are associated with penalty costs. To investigate the influence of the time windows on the planning results, we conduct experiments based on five different time windows, as depicted in **Figure 4**.

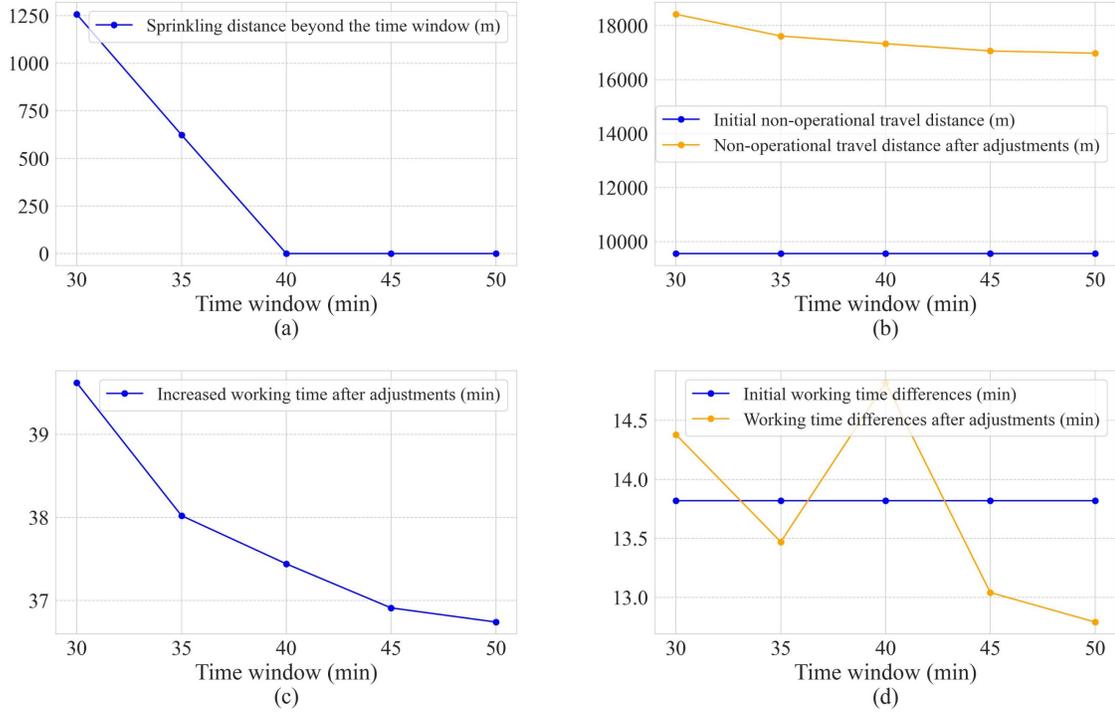

**Figure 4. Results under different time windows.**

**Figure 4**(a) shows that although sprinkling beyond the time window incurs significant penalty costs, when the time window duration is less than 40 minutes, there will still be some demand links that are difficult to be sprinkled within the time window. That is, for the instance 4, the 40-minute time window emerges as the critical threshold for ensuring the sprinkling of all real-time demand links when they are received 60 minutes after the start of the operation.

In terms of the impact of the time window duration on non-operational travel distance, as depicted in **Figure 4**(b). As the time window restrictions are relaxed, the non-operational travel distance after route adjustments gradually decreases, and this declining trend eventually stabilizes. Numerically, the influence of the time window on the non-operational travel distance is generally marginal. It can be inferred that when the time window exceeds 50 minutes, the impact of the time window restrictions on the planning results will be minimal, or even negligible. When only the time window duration changes, the total sprinkling distance remains constant. Consequently, the increased working time after route adjustments is determined by the increased non-operational travel distance. Therefore, the trend of the increased working time after route adjustments in **Figure 4**(c) parallels the trend of non-operational travel distance after route adjustments in **Figure 4**(b) as the time window duration varies.

For the working time differences among the sprinklers at different time windows, except for the 40-minute time window, the working time difference of the sprinklers diminishes as the time window duration increases. This may be because the 40-minute time window is the critical threshold. Furthermore, after route adjustments, the working time difference among the sprinklers could even be smaller than it was in the fixed routing plans. Therefore, extending the response time window proves advantageous for achieving a more balanced distribution of working time among the sprinklers.

## 5.4 Real-world deployment of the Dynamix-CARP solver

The proposed Dynamix-CARP solver is implemented in a decision support system named *Alpha MAPS* and deployed in an air quality management pilot project in a 3640 km$^2$ area in Chengdu, China. This project features real-time and responsive measures, enabled by pervasive air quality sensing (Dai and Han, 2023; Ji et al., 2023) and multi-source big data platforms (Song et al., 2020).

**Figure 5** illustrates how the proposed algorithm works with historical and real-time sprinkling demands within the system. In the sprinkler operation module of Alpha MAPS, the fixed sprinkling demands for the road segments are generated according to historical (static) information including road dust level, average traffic loads, and population density. Then, the Dynamix-CARP solver generates a fixed routing plan for all the sprinklers. On the other hand, the real-time (dynamic) demands are generated as follows. First, the system is instrumented with FRD sensing and evaluation capabilities, building on real-time traffic information, roadside FRD sensors, as well as vehicle patrols. When an instance of FRD pollution is detected, an inspector will be dispatched to confirm such FRD pollution on site, while reporting its coordinates and relevant information using an App specifically designed for on-site inspection. Once such a real-time demand is received, the system calls the Dynamix-CARP solver to make real-time adjustment of the current routes.

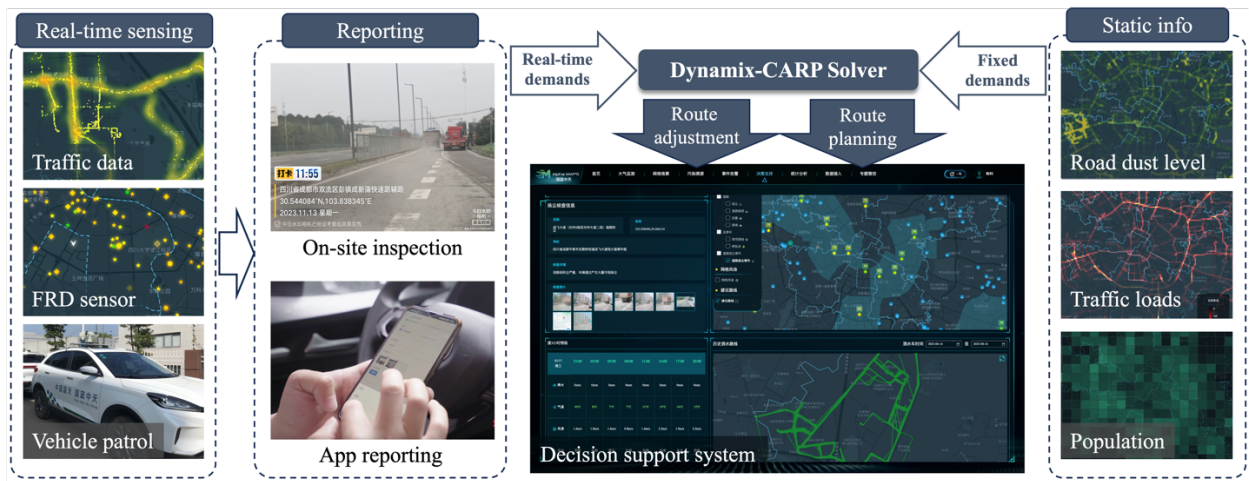

**Figure 5. Workflow of the Dynamix-CARP solver embedded in the Alpha MAPS system.**

## 6. Conclusions

This paper focuses on the routes planning of sprinklers on urban streets, considering the real-time new demands during sprinkling operations, and models it as a multi-depot mixed capacitated arc routing problem with real-time demands. Due to its complexity, an improved ALNS algorithm is introduced, which utilizes a tabu-list to prevent redundant searches and proposes an acceleration mechanism to reduce unnecessary search domains so that enhance the algorithm efficiency. Perturbations are also used to improve solutions quality and avoid local optima. Test results on instances demonstrate that the improved algorithm outperforms the ALNS algorithm in both computation time and solution quality, offering the potential for efficient and intelligent urban sprinkling planning.

Sensitivity analysis reveals that real-time route adjustments significantly increase the non-operational travel distance of the sprinklers, and the later the route adjustments occur, the greater this increase may be. Early or delayed route adjustments, as well as relaxing the sprinkling time window of the new demands, can reduce the working time difference among the sprinklers. But the time window has a comparatively minor impact on the route planning results.

In future research, it would be reasonable to explore and compare the effectiveness of other heuristic algorithms in solving this problem. Additionally, the model and algorithm could be tested with areas of different shapes to assess their adaptability and performance in diverse urban environments. Furthermore, investigating routes planning of sprinklers with multiple cycles in a day and real-time route adjustments to multiple scenarios would be beneficial in enhancing the sprinkling scheme's applicability and versatility.

**Declaration of generative AI and AI-assisted technologies in the writing process**

During the preparation of this work, the authors used ChatGPT in order to polish the writing of this paper. After using this tool, the authors reviewed and edited the content as needed and take full responsibility for the content of the publication.